\documentstyle[12pt,amssymb]{amsart}
\newcommand{\Ga}{\Gamma}
\newcommand{\th}{\theta}
\newcommand{\G}{{\Bbb G}}
\newcommand{\Tr}{\operatorname{Tr}}
\newcommand{\GL}{\operatorname{GL}}
\newcommand{\GSp}{\operatorname{GSp}}
\renewcommand{\Sp}{\operatorname{Sp}}
\newcommand{\coker}{\operatorname{coker}}
\newcommand{\ad}{\operatorname{ad}}
\newcommand{\PP}{{\cal P}}
\renewcommand{\gg}{{\frak g}}
\newcommand{\gl}{{\frak gl}}
\newcommand{\Mat}{\operatorname{Mat}}

\newcommand{\und}{\underline}
\newtheorem{thm}{Theorem}[section]
\newtheorem{prop}[thm]{Proposition}
\newtheorem{lem}[thm]{Lemma}

\theoremstyle{definition}

\numberwithin{equation}{section}
         
\newcommand{\Pf}{\noindent {\it Proof}}

\newcommand{\rk}{\operatorname{rk}}

\newcommand{\SL}{\operatorname{SL}}
\renewcommand{\a}{\alpha}
\renewcommand{\b}{\beta}
\newcommand{\Cone}{\operatorname{Cone}}
\newcommand{\si}{\sigma}
\newcommand{\im}{\operatorname{im}}
\newcommand{\ed}{\qed\vspace{3mm}}
\newcommand{\Pic}{\operatorname{Pic}}

\newcommand{\lan}{\langle}
\newcommand{\ran}{\rangle}
\renewcommand{\AA}{{\cal A}}
\newcommand{\FF}{{\cal F}}
\newcommand{\C}{{\Bbb C}}
\newcommand{\Z}{{\Bbb Z}}

\newcommand{\Lie}{\operatorname{Lie}}

\newcommand{\om}{\omega}
\newcommand{\ra}{\rightarrow}

\newcommand{\la}{\lambda}
\newcommand{\id}{\operatorname{id}}

\newcommand{\ot}{\otimes}

\renewcommand{\O}{{\cal O}}
\renewcommand{\P}{{\Bbb P}}

\newcommand{\Hom}{\operatorname{Hom}}
\newcommand{\End}{\operatorname{End}}
\newcommand{\Ext}{\operatorname{Ext}}

\newcommand{\MM}{{\cal M}}

\newcommand{\D}{{\cal D}}
\newcommand{\wt}{\widetilde}
\newcommand{\de}{\delta}

\newcommand{\sub}{\subset}
\newcommand{\lrar}[1]{\begin{picture}(50,10)(-25,-5)
\put(-25,0){\vector(1,0){50}}
\put(0,5){\makebox(0,0)[b]{\mbox{$#1$}}}
\end{picture}}

\newcommand{\ldar}[1]{\begin{picture}(10,50)(-5,-25)
\put(0,25){\vector(0,-1){50}}
\put(5,0){\mbox{$#1$}}
\end{picture}}

\title{Poisson structures and
birational morphisms associated with bundles on elliptic curves}
\author{A. Polishchuk}
         
\begin{document}
\maketitle
                   
\vspace{1cm}
   
Let  $X$ be a complex elliptic curve. 
In  this paper we define a natural Poisson structure on the
moduli  spaces  of  stable  triples  $(E_1,E_2,\Phi)$  where
$E_1,E_2$ are algebraic vector bundles on $X$ of fixed ranks
$(r_1,r_2)$ and degrees $(d_1,d_2)$,
$\Phi:E_2\ra E_1$ is a homomorphism. Such triples are considered
up to an isomorphism and the stability condition depends on
a real parameter $\tau$.
These moduli spaces were introduced by  S.~Bradlow  and
O.~Garcia-Prada \cite{BG}.
Our  Poisson structure induces a Poisson structure on similar
moduli spaces with fixed determinants of $E_1$ and $E_2$.
For    $E_2=\O$    and    some    values    of    parameters
$(r_1,r_2,d_1,d_2,\tau)$ the latter moduli
spaces are just the projective spaces. In particular,
one of these moduli spaces is $\P\Ext^1(E,\O)$, where
$E$ is a fixed stable bundle.
The  corresponding  Poisson  structures  on $\P\Ext^1(E,\O)$
were first defined and studied by B.~Feigin and A.~Odesskii \cite{FO}.
Moreover, they constructed a flat family of
quadratic algebras ({\it Sklyanin algebras})
$Q_{d,r}(x)$ where $d=\deg E$, $r=\rk E$, parametrized by $x\in X$
such that $Q_{d,r}(0)$ is the symmetric algebra in $d$
variables  and the quadratic Poisson bracket on the symmetric
algebra associated with this deformation induces the above
Poisson structure on $\P^{d-1}$. The algebra $Q_{d,r}(x)$ is
defined  as  the  associative  algebra  over $\C$
with $d$ generators $t_i, i\in\Z/d\Z$ and defining relations
\begin{equation}\label{relation}
\sum_{n\in\Z/d\Z}\frac{\th_{j-i+(r-1)n}(0)}
{\th_{j-i-n}(-x)\th_{rn}(x)} t_{r(j-n)}t_{r(i+n)}=0
\end{equation}
for $i,j\in\Z/d\Z$, where
$\th_m$, $m\in\Z/d\Z$ are certain 
theta-functions of level $d$ on $X$ (see \cite{FOf}).

For some other values of parameters we  get  as  moduli
space the projective space $\P H^0(E)$ where $E$ is
a stable bundle on $X$.
When the parameter $\tau$ changes continuosly
the corresponding moduli space doesn't change except when
$\tau$  passes  a finite number of rational values, in which
case the moduli space undergoes 
birational transformations (flips) compatible with Poisson structures.
In  particular,  if  we  start  with  a stable bundle $E$ of
rank $r$ and degree $d>r$ such that $d$ is relatively prime to
$r+1$, then we get a sequence of Poisson birational morphisms connecting
$\P\Ext^1(E,\O)$ and $\P H^0(E')$ where $E'$ is the unique
stable bundle of rank $r+1$ with $\det E'\simeq\det E$.
         
On the other hand, using Fourier-Mukai transform one constructs
an action of a central extension of $\SL_2(\Z)$ by $\Z$ on $D^b(X)$,
the derived category of coherent sheaves on $X$ (see \cite{Mukai}).
Using this action we construct for every stable bundle $E$ of rank $r$
and degree $d$ an
isomorphism (compatible with Poisson structures)
between the spaces $\P H^0(E)$ and
$\P\Ext^1(E',\O)$ where $E'$ is  certain  stable  bundle  of
degree $d$ and rank $r'$ satisfying the congruence relation
$r\cdot r'\equiv -1\mod(d)$.
 
Another application of $\SL_2(\Z)$-action gives a Poisson isomorphism
between $\P\Ext^1(E,\O)$ and $\P\Ext^1(E',\O)$ where
$E$ and $E'$ are stable bundles of the same degree $d$ whose ranks
$r$ and $r'$ satisfy the congruence $r\cdot r'\equiv 1\mod(d)$.
This reflects the fact noticed in \cite{FOf}
that the corresponding Sklyanin algebras 
$Q_{d,r}(x)$ and  $Q_{d,r'}(x)$ are isomorphic for any $x\in X$.
 
It turns out that
the above birational and regular isomorphisms of projective spaces
fit together in the following way. For every $d>0$ we can
consider the disjoint union of $(d-1)$-dimensional projective
spaces indexed by the set of residues $R_d\sub\Z/d\Z$ consisting
of $r$ such
that both $r$ and $r+1$ are relatively prime to $d$.
Namely, for every $r$ ($0<r<d$) the corresponding projective space
is $\P\Ext^1(E,\O)$ where $E$ is stable bundle of rank $r$ and
degree $d$. Then the above birational and regular isomorphisms
of projective spaces generate the birational action of $S_3$ (the group
of permutations in 3 letters) on this disjoint union. 
The corresponding action of $S_3$ on the set of connected components
$R_d$ is generated by operators $r\mapsto r^{-1}$ and
$r\mapsto -r-1$.

Finally, we show how to generalize our Poisson brackets to the
similar moduli stacks for other reductive groups.
Namely, we fix the following data: a reductive group $G$,
its finite-dimensional representation $V$ and a symmetric
$\gg$-invariant tensor $t\in (S^2\gg)^{\gg}$ where
$\gg$ is the Lie algebra of $G$. These data should satisfy
the following condition: the operator $t_*:S^2 V\ra S^2V$
induced by $t$ should be zero. Then we consider the
moduli stack of pairs $(P,s)$ where $P$ is a principal 
$G$-bundle on $X$, $s\in V(P)$ is a section of the associated
vector bundle. Given a trivialization of $\om_X$ we construct
a canonical Poisson structure on this moduli stack.
In the case when $G=\GL_{r_1}\times\GL_{r_2}$, $V$ is the
space of $r_1\times r_2$-matrices, there is a natural choice
of the tensor $t$ such that the above condition is satisfied
and we recover our Poisson structure on moduli of triples. 
The simplest case involving other groups than $\GL$ is
the case $G=\GSp(V)$ where $V$ is the symplectic vector
space, $\GSp(V)$ is the group of automorphisms preserving
the symplectic form up to a non-zero constant.       
In this case there is a canonical choice of $t$, so we get
a Poisson structure on the moduli stack of pairs
$(E,s)$ where $E$ is a vector bundle equipped with a symplectic
form $E\times E\ra L$ (where $L$ is a line bundle),
$s$ is a section of $E$. 
  
\section{Stable triples}
         
Let us recall the definition of stable triples
from \cite{BG}. Let $T=(E_1,E_2,\Phi)$ be a triple consisting
of two vector bundles $E_1$ and $E_2$ on $X$ and a homomorphism
$\Phi:E_2\ra E_1$. For a real parameter $\si$ the $\si$-degree
of $T$ is defined as follows:
$$\deg_{\si}(T)=\deg(E_1)+\deg(E_2)+\si\cdot\rk(E_2).$$
Now the $\si$-slope of $T$ is defined by the formula
$$\mu_{\si}(T)=\frac{\deg_{\si}(T)}{\rk(E_1)+\rk(E_2)}.$$
Note that if $L$ is a line bundle then we can define a tensor
of a triple $T$ with $L$ naturally, so that one has
$\mu_{\si}(T\otimes L)=\mu_{\si}(T)+\deg L$.
       
The triple $T$ is called $\si$-stable if for every
non-zero proper subtriple $T'\sub T$ one has
$\mu_{\si}(T')<\mu_{\si}(T)$.
Sometimes it is convinient to introduce another stability parameter
$\tau=\mu_{\si}(T)$.
         
The category of triples $T=(E_1,E_2,\Phi)$ is equivalent to the
category of extensions
\begin{equation}\label{equivext}
0\ra p^*E_1\ra F\ra p^*E_2(2)\ra0
\end{equation}
on $X\times\P^1$ where $p:X\times\P^1\ra X$ is the projection.
Indeed, the space of such extensions is
$\Ext^1_{X\times\P^1}(p^*E_2(2),p^*E_1)\simeq\Hom_X(E_2,E_1)$.
This extension has a unique $\SL_2$-equivariant structure
and as shown in \cite{BG} the $\si$-stability condition on $T$
is equivalent to the $\SL_2$-equivariant stability
of $F$ with respect to some polarization on $X\times\P^1$
depending on $\si$.
         
Let us denote by $\MM_{\si}=\MM_{\si}(d_1,d_2,r_1,r_2)$ the moduli space of
$\si$-stable triples $T=(E_1,E_2,\Phi)$ on $X$ with $\deg(E_i)=d_i$,
$\rk  E_i=r_i$.
When using another stability parameter $\tau$ we will denote
the same moduli space by $\MM_{\tau}$.
This  moduli space can be constructed using
geometric  invariant  theory  as  in  \cite{Be}.
       
We claim that in the case of elliptic  curve  all  these  moduli
spaces are smooth.
         
\begin{lem} The moduli space $\MM_{\si}$ is smooth.
\end{lem}
         
\Pf. According to \cite{BG} we have to show that
$H^2(X\times\P^1,\und{\End} F)^{\SL_2}=0$ for the
$\SL_2$-equivariant vector bundle $F$ associated with a
$\si$-stable triple. Consider the exact sequence
$$0\ra K\ra \und{\End} F\ra\und{\Hom}(p^*E_1,p^*E_2(2))\ra 0.$$
Since the direct image $Rp_*$ of the last term will have
no $\SL_2$-invariant part we have
$H^*(X\times\P^1,\und{\End} F)^{\SL_2}\simeq H^*(X\times\P^1,K)$.
Now $K$ sits in the exact triangle
$$\und{\Hom}(p^*E_2(2),p^*E_1))\ra K\ra\und{\End} p^*E_1\oplus
\und{\End} p^*E_2\ra\und{\Hom}(p^*E_2(2),p^*E_1)[1]\ra\ldots$$
It follows that equivariant direct image of $K$ with respect
to the projection $p$ is quasi-isomorphic to the complex
\begin{equation}\label{tangcomp}
C^{\cdot}:\und{\End} E_1\oplus\und{\End} E_2\stackrel{d}{\ra}
\und{\Hom}(E_2,E_1)
\end{equation}
concentrated in degrees $0$ and $1$, where $d(A,B)=A\Phi-\Phi B$.
We have the exact sequence of cohomologies
$$H^1(X,\und{\End}E_1\oplus\und{\End} E_2)\ra
H^1(X,\und{\Hom}(E_2,E_1))\ra H^2(X, C^{\cdot})\ra 0.$$
Now by Serre duality we have
$H^1(X,\und{\Hom}(E_2,E_1))^*\simeq H^0(X,\und{\Hom}(E_1,E_2))$.
According to Lemma 4.4 of \cite{BG} this space is zero
unless  $\Phi$  is  an  isomorphism.  In the latter case the
first arrow in the above exact sequence is surjective, so in
either case we get $H^2(X\times\P^1,K)=H^2(X, C^{\cdot})=0$.
\ed
         
The proof of this lemma also shows that the tangent space to
$\MM_{\si}$  at  a  triple  $T$  is  identified   with   the
hypercohomology space $H^1(X,C^{\cdot})$ where $C^{\cdot}$
is the complex (\ref{tangcomp}). This can also be
shown directly considering infinitesemal deformations
of the first order for triples.
         
\section{Poisson structure}
\label{mainsec}
      
Let   us  fix  a  trivialization  $\om_X\simeq\O_X$  of  the
canonical bundle of $X$. Then we can
define a Poisson structure on the moduli space of
triples $\MM_{\si}$. As we have seen above the tangent space
to $\MM_{\si}$ at a triple $T$ is identified with
$H^1(X,C)$ where $C$ is the complex
(\ref{tangcomp}). By Serre duality the cotangent space
is isomorphic to $H^0(X,C^*)=H^1(X,C^*[-1])$,
where the complex $C^*[-1]=((C^1)^*\stackrel{-d^*}{\ra}(C^0)^*)$ is
concentrated in degrees $0$ and $1$. Using the natural
autoduality of $\End E_i$ the complex $C^*[-1]$ can be identified with
$$\und{\Hom}(E_1,E_2)\stackrel{-d^*}\ra
\und{\End}E_1\oplus\und{\End}E_2$$
where $-d^*(\Psi)=(-\Phi\Psi,\Psi\Phi)$.
Now let us consider the
morphism of complexes $\phi:C^*[-1]\ra C$ with components
$\phi_1=0$ and
\begin{equation}\label{phi0}
\phi_0:\und{\Hom}(E_1,E_2)\ra
\und{\End}E_1\oplus\und{\End}E_2:
\Psi\mapsto (\Phi\Psi,\Psi\Phi).
\end{equation}
Since $d\circ\phi_0=0$,  we  have  indeed  the  morphism  of
complexes.  Therefore,  we  can  take  the  induced  map  on
hypercohomologies
$$H_T=\phi_*:H^1(X,C^*[-1])\ra H^1(X,C).$$
Note that we get a map from the cotangent space
to the tangent space of $\MM_{\si}$ at $T$.
This construction easily globalizes to give a morphism
$H$ from the cotangent  bundle  to  the  tangent  bundle  of
$\MM_{\si}$.
         
\begin{thm} $H$ defines a Poisson structure on $\MM_{\si}$.
\end{thm}
         
\Pf . Let us check that $H^*=-H$. First of all we claim that
$\phi^*[-1]=\phi$ in  the  homotopy  category  of  complexes.
Indeed, by definition $\phi^*[-1]$ has components
$(\phi^*[-1])_0=0$ and
$$(\phi^*[-1])_1:\und{\End}E_1\oplus\und{\End}E_2\ra
\und{\Hom}(E_2,E_1):(A,B)\mapsto A\Phi+\Phi B.$$
Now let us consider the map
$$h:\und{\End}E_1\oplus\und{\End}E_2\ra\und{\End}E_1\oplus
\und{\End}E_2:(A,B)\mapsto (-A,B).$$
Immediate check shows that $h$ provides a homotopy from
$\phi^*[-1]$ to $\phi$. Now the skew-commutativity of $H$
follows immediately from the skew-commutativity of the
natural pairing
$H^1(C)\otimes H^1(C')\ra H^2(C\otimes C')$
that comes from the minus sign in the commutativity
constraint for the tensor product of complexes.

The Jacoby identity will be proven in section \ref{gener}
in more general situation.
\ed
        
We can interpret the Poisson bivector $H$ in terms of
$\SL_2$-equivariant bundles on $X\times\P^1$ as follows.
Let  $F$  be the extension (\ref{equivext}) associated with a triple
$T$. Then the tangent space to $\MM_{\si}$ at $T$ is
identified with $H^1(X\times\P^1, \und{\End} F)^{\SL_2}$,
hence by Serre duality the cotangent space is identified with
$H^1(X\times\P^1, \und{\End}(F)(-2))$. Now we claim that $H_T$
is induced by some canonical morphism
$$\a:\und{\End}F(-2)\ra\und{\End}F$$
on $X\times\P^1$.
Namely, let $\und{\End}(F,p^*E_1)$  be  the  kernel  of  the
natural projection
$\und{\End}F\ra\und{\Hom}(p^*E_1,p^*E_2(2))$.
The dual morphism to the embedding gives a morphism
$\und{\End}F\ra\und{\End}(F,p^*E_1)^*$.
Thus, to construct $\a$ it is sufficient to construct a morphism
$$\wt{\a}:\und{\End}(F,p^*E_1)^*(-2)\ra\und{\End}(F,p^*E_1).$$
Now the bundle  $\und{\End}(F,p^*E_1)$ sits in the following
exact triple
$$0\ra\und{\Hom}(p^*E_2(2),p^*E_1)\ra \und{\End}(F,p^*E_1)
\ra\und{\End}(p^*E_1)\oplus\und{\End}(p^*E_2)\ra 0.$$
We have the morphism to the last term of this triple
$$p^*\phi_0:\und{\Hom}(p^*E_1,p^*E_2)\ra\und{\End}(p^*E_1)\oplus
\und{\End}(p^*E_2)$$
where $\phi_0$ is defined in (\ref{phi0}). It is easy to check
that $p^*\phi_0$ lifts uniquely to a morphism
$$\und{\Hom}(p^*E_1,p^*E_2)\ra \und{\End}(F,p^*E_1).$$
Now we define $\wt{\a}$ to be the composition of the latter
morphism with the natural projection
$\und{\End}(F,p^*E_1)^*(-2)\ra\und{\Hom}(p^*E_1,p^*E_2)$.
        
One has the natural morphism $\det:\MM_{\si}\ra\Pic(X)^2$
associating to a triple $(E_1,E_2,\Phi)$ the pair of line
bundles $(\det E_1,\det E_2)$. We claim that $\det$ is
a Casimir morphism, i.e. preimage of any local function
downstairs is a Casimir function upstairs (that is a function
having zero Poisson bracket with any other function).
Indeed, the cotangent map to $\det$ is just the natural
embedding
$$i:H^0(X,\O_X)^2\ra H^1(X, C^{\cdot})$$
which factors through $H^0(X,C^1)=H^0(X,\und{\End} E_1)\oplus
H^0(X,\und{\End} E_2)$. On the other hand, $H$ factors through
the map $H^1(X, C^{\cdot})\ra H^1(X,C^0)$, hence the image
of $i$ is killed by $H$.
In particular, the Poisson bracket on $\MM_{\si}$
induces Poisson brackets on the fibers of the morphism $\det$.
These fibers can be identified with moduli spaces
$\MM_{\si}(L_1,L_2,r_1,r_2)$ of triples with
fixed determinants $\det E_i\simeq L_i$.
Tensoring with a fixed line bundle $L$ gives a Poisson isomorphism
$\MM_{\si}(L_1,L_2,r_1,r_2)\simeq
\MM_{\si}(L_1\ot L^{\ot r_1}, L_2\ot L^{\ot r_2},r_1,r_2)$.
An automorphism $\phi:X\ra X$ induces an isomorphism of moduli
spaces $\MM_{\si}(L_1,L_2,r_1,r_2)\ra\MM_{\si}(\phi^*L_1,
\phi^*L_2,r_1,r_2)$ compatible with Poisson structures.
    
\section{Special cases}
        
For any bundle $E_1$ and a subbundle $E_2\sub  E_1$  let  us
denote   by   $\und{\End}(E_1,E_2)$   the   sheaf  of  local
homomorphisms of $E_1$ preserving $E_2$. In other words,
this is the kernel of the natural projection
$\und{\End}(E_1)\ra\und{\Hom}(E_2,E_1/E_2)$.
        
\begin{lem}\label{subbun}
Let $T=(E_1,E_2,\Phi)$
be a $\si$-stable triple such that $\Phi:E_2\ra E_1$ is an embedding of
$E_2$ as a subbundle. Then
the tangent space to $\MM_{\si}$ at $T$ can be identified
with $H^1(X,\und{\End}(E_1,E_2))$.
\end{lem}
                   
\Pf. The  natural embedding
$\und{\End}(E_1,E_2)\ra\und{\End}(E_1)\oplus\und{\End}(E_2)$
induces    the   quasi-isomorphism   $\und{\End}(E_1,E_2)\ra
C^{\cdot}$. Hence the assertion.
\ed
       
Under the identification of this lemma our Poisson structure
at the triple $T$ for which $\Phi$  is  an  embedding  of  a
subbundle can be described as follows. From the exact triple
$$0\ra\und{\Hom}(E_1/E_2,E_2)\ra\und{\End}E_1
\ra\und{\End}(E_1,E_2)^*\ra 0$$
we get a boundary homomorphism
$$H^0(X,\und{\End}(E_1,E_2)^*)\ra
H^1(X,\und{\Hom}(E_1/E_2,E_2)).$$
Composing it with the natural morphism
$$H^1(X,\und{\Hom}(E_1/E_2,E_2))\ra
H^1(X,\und{\End}(E_1,E_2))$$
we get a morphism
$$H^1(X,\und{\End}(E_1,E_2))^*\simeq
H^0(X,\und{\End}(E_1,E_2)^*)\ra H^1(X,\und{\End}(E_1,E_2)),$$
which  coincides  with  $H_T$  under  identification  of Lemma
\ref{subbun}. Equivalently, we may start with
the natural morphism
$$H^0(X,\und{\End}(E_1,E_2)^*)\ra
H^0(X,\und{\Hom}(E_2,E_1/E_2))$$ and compose it with the
boundary homomorphism
$$H^0(X,\und{\Hom}(E_2,E_1/E_2))\ra
H^1(X,\und{\End}(E_1,E_2))$$
coming from the exact triple
\begin{equation}\label{imptriple}
0\ra\und{\End}(E_1,E_2)\ra\und{\End} E_1\ra
\und{\Hom}(E_2, E_1/E_2)\ra 0.
\end{equation}
The  equivalence  of  this  description with the previous one
follows immediately from the commutative diagram
\begin{equation}
\begin{array}{ccccc}
\und{\Hom}(E_1/E_2,E_2) &\lrar{} &\und{\End}(E_1,E_2)
&\lrar{} &\und{\End} E_2\oplus\und{\End} E_1/E_2 \\
\ldar{\id} & & \ldar{} & & \ldar{} \\
\und{\Hom}(E_1/E_2,E_2) &\lrar{} &\und{\End} E_1
&\lrar{} &\und{\End}(E_1,E_2)^* \\
&&\ldar{} & & \ldar{}\\
&&\und{\Hom}(E_2,E_1/E_2) &\lrar{\id}&
\und{\Hom}(E_2,E_1/E_2)
\end{array}
\end{equation}
     
It is easy to check using the above descriptions
that the restriction of our Poisson bracket to the space of triples
for  which  $\Phi$  is  an  embedding  of  a  subbundle is a
particular case of the Poisson bracket on moduli of
principal bundles over
parabolic subgroups  defined  by  Feigin  and  Odesskii  in
\cite{FO}.
      
\begin{lem}\label{simple}
Keep the assumptions of Lemma \ref{subbun}.
Then we have an exact triple
$$0\ra \End E_2\oplus\End (E_1/E_2)\ra\ker H_T\ra
\End E_1/\End(E_1,E_2)\ra 0.$$
If $\Phi(E_2)$ is a direct summand of $E_1$ then $H$ vanishes
at $T$. Otherwise, the dimension of $\ker H_T$ is minimal
(and equals to $2$) if and only if the bundles $E_1$, $E_2$
and $E_1/E_2$ are simple.
\end{lem}
     
\Pf . Considering the second of the above descriptions of
$H_T$ we see immediately that there is an exact sequence
\begin{equation}\label{triplepf}
0\ra \End E_2\oplus \End (E_1/E_2)\ra\ker H_T\ra
\ker (H^0(X,\und{\Hom}(E_2,E_1/E_2))\ra
H^1(\und{\End}(E_1,E_2)).
\end{equation}
Using (\ref{imptriple}) the last term can be identified with
$\End   E_1/  \End(E_1,E_2)$. Moreover, the  last  arrow in
(\ref{triplepf}) is
surjective since we have a natural map
$\End E_1\ra\ker H_T$ coming from the morphism
$\und{\End} E_1\ra\und{\End}(E_1,E_2)^*$, and the composition
of this map with the last arrow of (\ref{triplepf})
is just the canonical projection to
$\End E_1/\End(E_1,E_2)$.
     
If $\Phi(E_2)$  is  a  direct  summand  in  $E_1$  then  the
boundary homomorphism used in the definition of $H_T$ is zero,
hence $H_T=0$. Otherwise, $\dim \ker H_T=2$ if and only if
$E_2$ and $E_1/E_2$ are simple and all endomorphisms of $E_1$
preserve $\Phi(E_2)$. We claim this can happen only when
$E_1$ is also simple. Indeed, let $A:E_1\ra E_1$ be any
endomorphism. By assumption $A$ preserves $\Phi(E_2)$.
Adding a constant to $A$ we may assume that
$A|_{\Phi(E_2)}=0$. Then it induces a map $E_1/E_2\ra E_2$.
However, $\sigma$-stability of our triple implies by Lemma 4.4
of \cite{BG} that $\Hom(E_1,E_2)=0$ since $\Phi$ is not an
isomorphism in our situation. It follows that
$\Hom(E_1/E_2,E_2)=0$,  hence,  $A=0$.  Thus, $\End E_1=\C$.
\ed
     
We are mainly interested in the case when $E_2=\O_X$, $\det E_1$
is fixed. In terms of parameter $\tau$ the stability condition
on $\Phi:\O_X\ra E_1$
is that $\mu(E'_1)<\tau$ for every proper non-zero subbundle
$E'_1\sub E_1$ and $\mu(E_1/E'_1)>\tau$ for every proper subbundle
$E'_1\sub E_1$ such that $\Phi\in H^0(X,E'_1)$.
        
Now let $E$ be a stable bundle on $X$ of degree $d$ and rank
$r$ (in particular, $gcd(d,r)=1$).
Set $\tau=\mu(E)=\frac{d}{r}$ and consider the moduli space
$\MM_{\tau}(\det E,\O_X,r+1,1)$. It is easy to see that the
stability condition on such a triple $\Phi:E_2=\O_X\ra E_1$ is
equivalent to the condition that $\Phi$ is nowhere vanishing
section and the quotient $E_1/\Phi(\O_X)$ is a stable bundle.
Moreover, since there exists a unique stable bundle of  rank
$r$ and determinant $\det E$ it follows that
$E_1/\Phi(\O_X)\simeq E$. Thus, we can identify the moduli
space of such triples with the projective space
$\P\Ext^1(E,\O_X)$.  If $gcd(d,r+1)=1$ then generic extension of  $E$  by
$\O_X$  is  stable.  Hence,  according to Lemma \ref{simple} in this case the
Poisson bracket on $\P\Ext^1(E,\O_X)$ is symplectic at general point.
 
Let  $t_x:X\ra X$ be the translation by
$x\in X$. Then by functoriality we have a natural Poisson
isomorphism
$$\P\Ext^1(E,\O_X)\wt{\ra} \P\Ext^1(t_x^*E,\O_X).$$
Note that as $x$ varies $t_x^*E$  runs  through  all  stable
bundles of given degree $d$ and rank $r$. Let $K\sub X$ be
the finite subgroup of order $d^2$ consisting of $x$ such that
$t_x^*\det E\simeq\det E$.
Then for $x\in K$ one has $t_x^*E\simeq E$, therefore, $K$
acts on $\P\Ext^1(E,\O_X)$ by linear transformations
preserving  the Poisson structure.
       
Another special moduli space associated with a fixed  stable
bundle $E$ is $\MM_{\tau}(\det E,\O_X,\rk E, 1)$ where
$\tau=\mu(E)$. Then the condition on a triple just means that
$E_1$ is stable, hence isomorphic to $E$, and $\Phi$ is
an arbitrary non-zero section. Therefore, this moduli space
can be identified with $\P H^0(X,E)$.
The Poisson bracket in this case can be described as follows.
The tangent space $T_{[s]}$ to the line generated by $s\in H^0(X,E)$ is
identified with $\coker(H^0(X,\O_X)\stackrel{s}{\ra}
H^0(X,E))$. Thus the cotangent space is
$T^*_{[s]}=\ker((H^1(X,E^*)\ra  H^1(X,\O_X))$.  Let  $D\sub  X$ be the
divisor of zeroes of $s$,  so  that  $s:\O_X\ra  E$  factors
as  $\O_X\ra\O_X(D)\ra  E$  where $\O_X(D)$ is embedded as a
subbundle into $E$. Then we have the natural map
\begin{equation}\label{map}
H^1(X,E^*)\ra H^1(X,E^*(D))\ra
H^1(X,\und{\End}(E,\O_X(D))).
\end{equation}
The latter space sits in the exact sequence
$$0\ra H^0(X,E(-D)/\O_X)\ra
H^1(X,\und{\End}(E,\O_X(D)))\ra H^1(X,\und{\End} E)
\simeq H^1(X,\O_X).$$
It follows that (\ref{map}) induces a map
$T^*_{[s]}\ra H^0(X,E(-D)/\O_X)$. Furthermore, it is easy to
check that its composition with the boundary homomorphism
$H^0(X,E(-D)/\O_X)\ra H^1(X,\O_X)$ is zero, hence, we get
a map
$$T^*_{[s]}\ra H^0(X,E(-D))/H^0(X,\O_X).$$
Now the latter space is naturally embedded into $T_{[s]}$
and the composition with this embedding gives our Poisson
structure at $[s]\in\P H^0(X,E)$.
     
\begin{lem} Let $H_{[s]}:T_{[s]}^*\ra T_{[s]}$ be the above
Poisson structure on $\P H^0(X,E)$. If $\rk E=1$ then
$H_{[s]}=0$. Otherwise, one has an exact sequence
$$0\ra H^1(X,\ad(E/\O_X(D)))\ra\coker H_{[s]}\ra
H^0(D,E|_D)\ra 0$$
where $D$ is the zero divisor of $s$.
\end{lem}
                      
\Pf . Let us denote by $V\sub T_{[s]}$ the subspace
$H^0(X,E(-D))/H^0(X,\O_X)$. Since the image of $H_{[s]}$
is contained in $V$ we have the exact sequence
$$0\ra\coker(T^*_{[s]}\ra V)\ra\coker H_{[s]}\ra
T_{[s]}/V\ra 0.$$
Since $E(-D)$ is stable of positive slope, it follows that
$H^1(X, E(-D))=0$, hence we have an isomorphism
$$T_{[s]}/V\simeq H^0(X,E)/H^0(X,E(-D))\simeq H^0(D,E|_D).$$
Now the assertion follows easily from the exact sequence
$$T_{[s]}^*\ra H^1(X,\und{\End}(E,\O(D)))\ra
H^1(X,\und{\End}(E/\O_X))\oplus H^1(X,\O_X)\ra 0.$$
\ed
 
More generally, we can consider the moduli space $\MM_{\tau}(L,\O_X,r,1)$
where $L$ is a fixed line bundle of degree $d$,
$\tau=\frac{d}{r}+\epsilon$ where $\epsilon>0$ is sufficiently small.
Then we get the moduli space of pairs $s:\O_X\ra E$
where $E$ is a semistable bundle with determinant $L$, $\rk E=r$, $s$ is
a section which doesn't belong to any destabilizing subbundle of $E$. We have
a Casimir morphism from this moduli space to the moduli stack of semistable
bundles, so the fibers inherit the Poisson structure. In particular,
if we take the semistable bundle $E=(E_0)^{\oplus k}$ where $E_0$ is a
stable bundle, then the corresponding fiber is identified with the
Grassmannian $G(k,H^0(E_0))$ of $k$-dimensional subspaces in $H^0(E_0)$,
so we get some family of Poisson structures on the Grassmannians.
     
\section{Fourier transforms}
                      
Let  $m:X\times X\ra X$ be the group law on $X$, $x_0\in X$ be
the neutral element.
Let
$$\PP=m^*\O_X(x_0)\otimes      p_1^*\O_X(-x_0)\otimes
p_2^*\O_X(-x_0)$$
be the Poincar\'e line bundle on $X\times X$
inducing an isomorphism of $X$ with the dual elliptic
curve. We denote by $\FF$ the corresponding Fourier-Mukai
transform which is an autoequivalence of the
the derived category $\D^b(X)$ of coherent sheaves on $X$ given
by
$$\FF(E)=Rp_{2*}(Lp_1^*E\sideset{^L}{}{\otimes}\PP).$$
One has $\FF\circ\FF\simeq (-\id_X)^*[-1]$ (see \cite{Mukai}).
It is easy to see that for every $E\in\D^b(X)$ one has
$\rk\FF(E)=\deg E$, $\deg\FF(E)=-\rk E$.
It follows that if $T:\D^b(X)\ra\D^b(X)$ is a composition of
some sequence of Fourier transforms and
tensorings with line bundles, then the vector
$v(T(E))=(\rk T(E), \deg T(E))$ is obtained from the vector
$v(E)=(\rk E, \deg E)$ by applying some matrix $A\in\SL_2(\Z)$.
Furthermore, one can lift the natural action of $\SL_2(\Z)$
on vectors $(\rk, \deg)$ to the action of a central extension
of $\SL_2(\Z)$ by $\Z$ on $\D^b(X)$.
More precisely, we can consider the
standard presentation of $\SL_2(\Z)$ by generators
$S=\left( \matrix 0 & 1\\ -1 & 0 \endmatrix \right)$ and
$R=\left( \matrix 1 & 0\\ 1 & 1 \endmatrix \right)$
subject to relations
$$S^2=(RS)^3,\ S^4=1.$$
Then the central extension in question is the group $\wt{\SL}_2(\Z)$
generated by $S$ and $R$ with the only relation $S^2=(RS)^3$.
The action of this group on $\D^b(X)$ is the following:
$S$ acts as the Fourier-Mukai transform while
$R$ acts as tensoring with $\O_X(x_0)$ (see \cite{Mukai}).
 
We will consider the action of $\wt{\SL}_2(\Z)$
on morphisms of stable bundles.
For this it will be useful to know the orbits of $\SL_2(\Z)$
on pairs of primitive vectors in $\Z^2$. First of all,
for a pair of vectors
$v_1=(r_1,d_1)$, $v_2=(r_2,d_2)$ such that $gcd(r_i,d_i)=1$ for $i=1,2$
the determinant $\det(v_1,v_2)\in\Z$ is invariant of $\SL_2(\Z)$.
We consider only pairs $v_1,v_2$ with $\det(v_1,v_2)\neq 0$.
For such pairs there is a second $\SL_2(\Z)$-invariant
$\a(v_1,v_2)\in(\Z/\det(v_1,v_2))^*$ defined from the condition
$$v_1\equiv\a(v_1,v_2)v_2\mod\det(v_1,v_2)\Z^2$$
It is easy to see that the $\SL_2(\Z)$-orbit of such $(v_1,v_2)$
consists of all pairs with the same $\det$ and $\a$.
                     
Henceforward, we restrict ourself to considering stable  bundles
on  $X$  with
determinant isomorphic to $\O_X(nx_0)$  for  some  $n$.  The
reason is that the group $\wt{\SL}_2(\Z)$
preserves the set $S_{x_0}$ of objects of the form
$E[k]$ where $k\in\Z$, $E$ is either a stable bundle
with determinant $\O_X(nx_0)$ for some $n$,
or $\O_{x_0}$. An element of $S_{x_0}$ is determined by
its degree and rank uniquely up to a shift. It follows that
the group $\wt{\SL}_2(\Z)$ acts transitively
on $S_{x_0}$. Furthermore, an element $T\in\wt{\SL}_2(\Z)$
is completely determined by its action on a pair of 
elements of $S_{x_0}$ which are not isomorphic up to shift.
  
The   first   immediate   consequence   of   the  action  of
$\wt{\SL}_2(\Z)$ is that for stable bundles $E_1$ and $E_2$
such that $\det E_1\simeq\det E_2\simeq\O_X(dx_0)$ and
$\rk  E_1\equiv  \rk  E_2\mod(d)$  there  is   a   canonical
isomorphism
$$\P\Ext^1(E_1,\O_X)\simeq\P\Ext^1(E_2,\O_X).$$
Indeed, under these conditions there is a unique element
$T\in\wt{\SL}_2(\Z)$ such that $T(\O_X)\simeq\O_X$ and
$T(E_1)\simeq E_2$. Considering the action of $T$
on morphisms from $E_1$ to $\O_X[1]$ we get the above
isomorphism.
  
Now for  every  stable  bundle  $E$  with $\det
E\simeq\O_X(dx_0)$ where $d>1$ we can find an element
$T\in\wt{\SL}_2(\Z)$ such that $T(E)\simeq \O_X$.
Then $T(\O_X)=E'[n]$ for
some stable bundle $E'$ and some $n\in\Z$. Since $\Hom(\O_X,E)\neq 0$
we should have $\Hom(T(\O_X), T(E))\neq 0$, hence $n=0$ or $-1$.
Consider first the case $n=-1$.
Then one has $\det E'\simeq \O_X(dx_0)$ and
$$r\cdot r'\equiv -1\mod(d)$$
where $r=\rk E$, $r'=\rk E'$
(this is deduced comparing invariants $\a$ for the pair of vectors
$(v(\O_X),v(E))$ and its image under $T$).
Conversely, for every $E'$ satisfying these conditions there
exists  a unique element $T\in\wt{\SL}_2(\Z)$ sending $E$ to
$\O_X$ and $\O_X$ to $E'[-1]$.
The transformation $T$ induces an isomorphism
\begin{equation}\label{linisom}
T_*:\P H^0(E)\wt{\ra}\P\Ext^1(E',\O_X).
\end{equation}
(an isomorphism of this kind with $r=1$, $r'=d-1$
was constructed in \cite{FMW} by a different method.)
Note that in the previous section we identified both sides of the
isomorphism (\ref{linisom})
with some special moduli spaces of pairs, in particular,
they carry natural Poisson structures.
     
\begin{prop} The isomorphism $T_*$ is compatible with Poisson structures.
\end{prop}
    
\Pf . Let   $s:\O_X\ra   E$   be   a   non-zero   section,
$\O_X\ra\wt{E}\ra E'$ be the corresponding extension of $E'$
by  $\O_X$ with class $T_*(s)\in\Ext^1(E',\O_X)$.
It  suffices to prove that $T_*$ preserves Poisson structure
over a non-empty open  subset,  hence  we  can  assume  that
$E/\O_X$ has no torsion.
    
By Serre duality the  cotangent  space
$T^*_{[s]}\P  H^0(E)$ can be identified with
$$\ker(\Ext^1(E,\O_X)\ra  H^1(\O_X))\simeq
\Ext^1(E/\O_X,\O_X)/\C\cdot e$$
where $e\in\Ext^1(E/\O_X, \O_X)$ is the class of the extension
$\O_X\ra E\ra E/\O_X$. Under this identification the Poisson bracket
on $\P H^0(E)$ at the point $[s]$ is induced by the natural morphism
$$H_{[s]}:\Ext^1(E/\O_X,\O_X)/\C\cdot e\ra T_{[s]}\P H^0(E)$$
which comes from the identification of $T_{[s]}\P H^0(E)$
with $\ker(H^1(\und{\End}(E,\O_X))\ra H^1(\O)^2)$ and
from the natural morphism
$$\Ext^1(E/\O_X,\O_X)=H^1(\und{\Hom}(E/\O_X,\O_X))\ra
H^1(\und{\End}(E,\O_X)).$$
In other words, we have a morphism from a neighborhood of $[e]$
in the space of extensions
$\P\Ext^1(E/\O_X,\O_X)$  to  $\P  H^0(E)$  (since   in   the
neighborhood of $[e]$ such an extension is necessarily
isomorphic to $E$), and $H_{[s]}$ is just the
tangent map to this morphism at the point $[e]$.
    
Similarly, the Poisson bracket on $\P\Ext^1(E',\O_X)$
at the point $[T_*s]$ can be identified with
the tangent map
$$\Hom(\wt{E},E')/\C\cdot f\ra
T_{[T_*(s)]}\P\Ext^1(E',\O_X)$$
to the local morphism from
$\P\Hom(\wt{E},E')$ to $\P\Ext^1(E',\O_X)$ at the point
$[f]$ where $f:\wt{E}\ra E'$ is the canonical morphism.
Here we use the natural identification of
$\Hom(\wt{E},E')/\C\simeq\ker(H^0(E')\ra H^1(\O_X))$
with the cotangent space to $\P\Ext^1(E',\O_X)$ at
$[T_*(s)]$.
Now  we  have  the  following  commutative square of local
morphisms in the neighborhood of points induced by $s$:
\begin{equation}
\begin{array}{ccc}
\P\Ext^1(E/\O_X,\O_X) & \lrar{} & \P H^0(E)\\
\ldar{T_*} && \ldar{T_*}\\
\P\Hom(\wt{E},E') &\lrar{} & \P\Ext^1(E',\O_X)
\end{array}
\end{equation}
Considering the corresponding commutative square of
tangent maps we get the compatibility of $T_*$ with
Poisson brackets.
\ed
  
For some other choice of autoequivalence  $T:\D^b(X)\ra\D^b(X)$
sending $E$ to $\O_X$ one has $T(\O_X)\simeq (E'')^*$ for a
stable  bundle  $E''$ of degree $\deg E''=\deg E=d$ and rank
$r''$ satisfying the congruence
$$r''\cdot r\equiv 1\mod(d)$$
where $r=\rk E$.
In this case we get an isomorphism
$$\P H^0(E)\simeq \P H^0(E'').$$
We claim that it is also compatible with the natural Poisson
structures on both sides. Indeed, this is proven exactly
as in the above proposition using the following commutative
diagram of local morphisms:
\begin{equation}
\begin{array}{ccc}
\P\Ext^1(E/\O_X,\O_X) & \lrar{} & \P H^0(E)\\
\ldar{T_*} && \ldar{T_*}\\
\P\Hom(E'',E''/\O_X) &\lrar{} & \P H^0(E'')
\end{array}
\end{equation}
    
Combining these isomorphisms we also get Poisson isomorphisms
between $\P\Ext^1(E_1,\O_X)$ and $\P\Ext^1(E_2,\O_X)$ for
stable bundles $E_1$ and $E_2$ of the same degree $d$ and
ranks $r_1$ and $r_2$ satisfying $r_1 r_2\equiv 1\mod (d)$.
    
We denote by $\MM(d,r)$ the projective space
$\P\Ext^1(E,\O_X)$  where  $E$  is  a  stable  bundle   with
determinant $\O_X(dx_0)$ and
rank $r$ (in particular, $gcd(d,r)=1$), considered as a Poisson
variety. Then the above results show
that $\MM(d,r)$ depends only on $d$ and on the residue
$r\in(\Z/d\Z)^*$, furthermore, we have
\begin{equation}\label{fourisom}
\MM(d,r)\simeq\MM(d,r^{-1}).
\end{equation}
Also for every stable bundle $E$ of degree $d>0$ and rank $r$
we have an isomorphism of Poisson varieties
$$\P H^0(E)\simeq\MM(d,-r^{-1}).$$
                                   
The  Poisson  isomorphism  (\ref{fourisom}) is the classical
limit of the following isomorphism of Sklyanin algebras
$$Q_{d,r}(x)\simeq Q_{d,r'}(x)$$
for every $\tau\in X$, where $rr'\equiv 1\mod(d)$ 
(cf. \cite{FOf}).
To see this isomorphism let us make the substitutions
$i=r'j'$, $j=r'i'$, $n=r'(n'+i'-j')$
in the
defining relation (\ref{relation}) of $Q_{d,r}(x)$.
Then using the relation $\th_{-i}(-x)=a\cdot b^i\cdot\th_i(x)$
(where $a$ and $b$ are some non-zero constants, $b^d=1$)
we can rewrite the quadratic relations in the form
\begin{equation}\label{relation2}
\sum_{n'\in\Z/d\Z}\frac{\th_{j'-i'-(r'-1)n'}(0)}
{\th_{r'n'}(x)\th_{j'-i'-n'}(-x)}
t_{j'-n'}t_{i'+n'}=0.
\end{equation}
Now it is obvious 
that the map $t_i\mapsto t_{r'i}$ defines an isomorphism
from $Q_{d,r}(x)$ to $Q_{d,r'}(x)$ as required.
   
\section{Birational transformations}
    
If one changes the stability parameter $\tau$ the moduli spaces
$\MM_{\tau}(L_1,L_2,r_1,r_2)$ undergo some birational transformations,
see \cite{T}, \cite{BDW}. Clearly,
these birational transformations are compatible with the Poisson
structures on their domain of definition.
In particular, considering moduli of pairs $\O_X\ra E$ where
$\deg E=d$, $\rk E=r$ are such that $gcd (r,d)=1$ and $gdc(r+1,d)=1$
we get a Poisson birational transformation from $\MM(d,r)$ to
$\MM(d,-(r+1)^{-1})\simeq\P H^0(E)$ where $E$ is a stable bundle
of degree $d$ and rank $r+1$.
Let us denote by $R_d\sub\Z/d\Z$ the set of residues $r$ such that
$r\in(\Z/d\Z)^*$ and $r+1\in(\Z/d\Z)^*$.
The map $\phi:r\mapsto -(r+1)^{-1}$ preserves $R_d$
and satisfies $\phi^3=\id$.
On the other hand, the involution
$\b:R_d\ra R_d: r\mapsto r^{-1}$ also preserves $R_d$ and
we have $\phi\b=\b\phi^{-1}$. It follows that $\b$ and $\phi$
generate the action of the symmetric group $S_3$ on $R_d$.
  
Recall that in the previous section we defined an isomorphism
$\MM(d,r)\ra\MM(d,r^{-1})=\MM(d,\b(r))$.
   
\begin{thm}\label{s3} The birational morphisms
$\MM(d,r)\ra\MM(d,\phi(r))$ and $\MM(d,r)\ra\MM(d,\b(r))$
defined above extend to an action of $S_3$ on
$\sqcup_{r\in R_d}\MM(d,r)$.
\end{thm}
   
\Pf . For every residue $r\in R_d$ let us denote by $E_r$ a
stable bundle with determinant $\O_X(dx_0)$ and rank
$\rk E\equiv r\mod(d)$ such that $0<\rk E<d$.
  
Let us check the relation $\phi^3=\id$. For this we have
to show that the corresponding composition of
birational transformations
$$\MM(d,r)\ra\MM(d,\phi(r))\ra\MM(d,\phi^2(r))\ra\MM(d,r)$$
is the identity. By definition the first arrow
is the composition of the map that associates to a
generic morphism $f:E_r\ra\O_X[1]$ the corresponding morphism
$\Cone(f)[-1]:\O_X\ra E_{r+1}$ (considered up to constant) with the
autoequivalence $T_r\in\wt{\SL}_2(\Z)$ which sends
$E_{r+1}$ to $\O_X[1]$ and $\O_X$ to $E_{\phi(r)}$.
It  follows  that  the  above  triple composition amounts to
applying the construction  $\Cone(\cdot)[-1]$  thrice  (note
that in our situation this construction is functorial) and
applying the functor $T_{\phi^2(r)}T_{\phi(r)}T_r$.
The triple composition of $\Cone(\cdot)[-1]$ is isomoprhic to
the shift $\id[-2]$, while
$T_{\phi^2(r)}T_{\phi(r)}T_r=\id[2]$, hence the assertion.
  
It remains to check the relation between our birational
transformations corresponding to the relation
$\phi\b=\b\phi^{-1}$. This amount to checking the following
relation between contravariant  functors  from  $D^b(X)$  to
itself:
$$DT_{\b\phi(r)}DU_{\phi(r)}T_r=U_r$$
where $D:\D^b(X)^{op}\ra\D^b(X)$ is the duality functor
$E\mapsto R\und{\Hom}(E,\O_X)$, $U_r\in\wt{\SL}_2(\Z)$
is the unique element sending $E_r$ to $\O_X$ and $\O_X$
to $E_{\b(r)}^*$. Note that for every $T\in\wt{\SL}_2(\Z)$
we have $DTD\in\wt{\SL}_2(\Z)$ (this follows from the
compatibility between the Fourier-Mukai transform and duality,
cf. \cite{Mukai}), moreover the corresponding involution on
$\SL_2(\Z)$ is just the conjugation by the matrix
$\left(\matrix 1 & 0 \\ 0 & -1 \endmatrix \right)$. Using this one can check
the above identity up to shift. Finally, since both parts send
$E_r$ to $\O_X$ the identity follows.
\ed
   
It follows from the above theorem that
for every $\sigma\in S_3$ and $r\in R_d$ such that $\si(r)=r$
we get a birational automorphism $f_{\sigma}$ of $\MM(d,r)$.
Since $\b$
acts as an isomorphism on our moduli spaces
there are essentially two different
cases to consider: the residue $r\in R_d$ is fixed $\phi$
or by $\phi\b$. The fixed points of $\phi$ are residues $r$ satisfying
$$r^2+r+1\equiv 0\mod(d).$$
In  this  case  we  get  a  Poisson  birational automorphism
$f_{\phi}$ of $\MM(d,r)$ such that $f_{\phi}^3=\id$.
The map $\phi\b$ has  the  only  fixed  point  $r=-2\mod(d)$
(provided that $d$ is odd), so we get a Poisson birational
involution $f_{\phi\b}$ of $\MM(d,d-2)$.
                                 
Let   us   describe   these  birational  automorphisms  more
explicitly.
Consider first the case when $r^2+r+1\equiv 0\mod(d)$, i.e.
$\phi(r)=r$.
Using the notation of the proof of the above theorem
we have $T_r(E_{r+1})=\O_X[1]$, $T_r(\O_X)=E_r$, and
$T_r(E_r)= E_{r+1}[1]$.
Now let us consider the closed subvariety of
$$Z\sub \P\Hom(E_{r+1},E_r)\times\P\Ext^1(E_r,\O_X)\times
\P\Hom(\O_X,E_{r+1})$$
consisting of triples of lines $([v_1],[v_2],[v_3])$
such that all three pairwise compositions $v_2\circ v_1$,
$v_3\circ v_2$ and $v_1\circ v_3$ are zeroes.
It is easy to see that any of three projections of $Z$ to
the projective spaces are birational. In fact, general point
of $Z$ corresponds to the exact triangle
$$\O_X\stackrel{v_3}{\ra} E_{r+1}\stackrel{v_1}{\ra} E_r
\stackrel{v_2}{\ra}\O_X[1].$$
In particular, $Z$ is birational to
$\MM(d,r)=\P\Ext^1(E_r,\O_X)$. On the other hand,
the functor $T_r$ gives rise to the automorphism of $Z$ given by
$$\Phi:(v_1,v_2,v_3)\mapsto (T_r(v_2)[-1],T_r(v_3),T_r(v_1)[-1]).$$
Now $\Phi$ induces our birational automorphism of $\MM(d,r)$ with cube
equal to the identity.
In the simplest non-trivial case $d=7$, $r=2$ the functor
we use has form $T_r:A\mapsto \FF(A(-2x_0))(3x_0)[1]$.
In  this  case one can describe $Z$ as the double blow-up of
$\MM(7,2)$. For this  it  is  more  convenient  to  use  the
isomorphism of $\MM(7,2)$ with $\P H^0(E_3)$.
Then we have a natural embedding
$S^2 C\ra\P H^0(E_3): D\mapsto H^0(E_3(-D))$ where $D$ is an effective divisor
of  degree  2 on $X$. We also define the 4-dimensional
variety $V\sub \P H^0(E_3)$ containing $S^2 C$ as
the union of chords of  $S^2  C\sub  \P  H^0(E_3)$  connecting
$D_1$   and   $D_2$   in   $S^2   C$   such   that  $D_1\cap
D_2\neq\emptyset$.
Then  our variety $Z$ is obtained by first blowing-up of $\P
H^0(E_3)$ along
$S^2 C$, and then blowing up the proper transform of $V$.
The automorphism $\Phi$ cyclically permutes the following three
divisors on $Z$: two exceptional  divisors  and  the  proper
transform of the chord
variety of $S^2 C$ (which is a hypersurface in $\P H^0(E_3)$).
  
In  the  case  $r\equiv  -2\mod(d)$  (where  $d$ is odd) the
birational autoequivalence of $\MM(d,d-2)$ is  described  as
follows.  Again  using  the  notation  from  the  proof   of
Theorem \ref{s3} we have
$U_{d-1}(E_{d-1})=\O_X$, $U_{d-1}(\O_X)=E_{d-1}^*$ and
$U_{d-1}(E_{d-2})=E_{d-2}^*[1]$.
Now let us consider the subvariety
$$Y\sub \P\Hom(E_{d-1},E_{d-2})\times\P\Ext^1(E_{d-2},\O_X)$$
consisting of pairs $([v_1],[v_2])$ such that
the composition $v_2v_1\in\Ext^1(E_{d-1},\O_X)$ is zero.
Then both  projections  of  $Y$  to  projective  spaces  are
birational.  On  the  other hand, the functor $U_{d-1}$ induces an
involution of $Y$ sending $(v_1,v_2)$ to
$(U_{d-1}(v_2)^*[1], U_{d-1}(v_1)^*[1])$,
hence our birational involution of
$\MM(d,d-2)=\P\Ext^1(E_{d-2},\O_X)$.
  
\section{Generalization}\label{gener}
      
In this section we consider a generalization of the main construction
of section \ref{mainsec} to the case of principal bundles with
other structural groups than $\GL$. For this note that the
datum of a triple $(E_1,E_2,\Phi)$ with  $\rk  E_i=r_i$  for
$i=1,2$ is equivalent to that of a pair $(P,s)$
where $P$ is a principal bundle with structure group
$\GL_{r_1}\times\GL_{r_2}$ and a section $s\in V(P)$
of the vector bundle $V(P)$ associated with the natural
representation of $\GL_{r_1}\times\GL_{r_2}$
on the space of $r_1\times r_2$-matrices.
      
To generalize this
let us consider a general reductive group $G$ and
its representation $V$. Then one can consider the moduli stack
$\MM_{G,V}$ of pairs $(P,s)$ where $P$ is a principal $G$-bundle,
$s\in H^0(X,V(P))$ is a global section of the corresponding
vector bundle associated to $V$ and $P$.
Let $\gg$ be the Lie algebra of $G$.
Assume that we are given a symmetric invariant tensor
$t\in S^2(\gg)^{\gg}$. Then $t$ induces a morphism of
$\gg$-modules $t_*:S^2(V)\ra S^2(V)$ as follows.
Let $t=\sum_i x_i\otimes y_i$, then
$$t_*(v\otimes v)=\sum (x_i\cdot v)\otimes (y_i\cdot v)$$
where $\cdot$ denotes the action of $\gg$ on $V$.
      
Assume that $t_*=0$. Then fixing a trivialization
$\om_X\simeq\O_X$ we can construct a Poisson bracket
on the smooth locus of $\MM$ as follows.
The tangent space to $\MM$ at a point $(P,s)$ can be identified with
the hypercohomology space $H^1(X,C^\cdot)$
where $C^{\cdot}$ is the complex
$\gg(P)\stackrel{d}{\ra} V(P)$ concentrated in degrees $0$ and $1$,
where $\gg(P)$ is the vector bundle associated with $P$ and the
adjoint representation, the map $d$ is induced by the
Lie action of $\gg(P)$ on $V(P)$: $d(A)=A\cdot s$.
Hence, the contangent space can be identified with
$H^1(C^*[-1])$. Now we can construct the morphism
of complexes $\phi:C^*[-1]\ra C$ as before
setting $\phi_1=0$ and
$\phi_0:V^*(P)\ra\gg(P)$ to be the composition of
$d^*$ with the map $\gg^*(P)\ra\gg(P)$ induced by $t$.
We claim that our condition on $t$ and $V$
implies that $d\circ\phi_0=0$ and that the obtained morphism
$H$ from the cotangent space of $\MM$ to the tangent space is
skew-symmetric. Indeed, essentially we have to check
that for every $v\in V$ the following composition is zero:
$$V^*\stackrel{d_v^*}{\ra}\gg^*\stackrel{t}{\ra}\gg
\stackrel{d_v}{\ra} V$$
where $d_v(A)=A\cdot v$. This is equivalent to the condition
$t_*(v\otimes v)=0$. Now the skew-symmetry follows as before:
the homotopy beteen $\phi$ and $\phi^*[-1]$ is constructed
using the map $\gg^*(P)\ra\gg(P)$ induced by $t$.

\begin{thm} The above construction defines
a Poisson bracket on the smooth locus of $\MM$.
\end{thm}

\Pf . We have to check the Jacoby identity for our bracket.
We  will  use  the  approach  similar  to that of \cite{Bo},
\cite{Bo2}. The Jacoby identity can be rewritten in terms of the morphism
$H:T_{\MM}^*\ra T_{\MM}$ as follows:
\begin{equation}\label{Jacobi}
H(\om_1)\cdot\lan H(\om_2),\om_3\ran-\lan [H(\om_1),H(\om_2)],
\om_3\ran+ cp(1,2,3)=0
\end{equation}
where $\om_i\in T_{\MM}^*$ are local 1-forms on $\MM$,
$[\cdot,\cdot]$ is the commutator of vector fields,
$cp(1,2,3)$ indicates terms obtained by cyclic permutation
of 1,2 and 3 from the first two terms.

Working over an affine \'etale open $U\ra\MM$ we can represent
every 1-form $\om\in T_{\MM}^*(U)$ by a Cech cocycle
$(\phi_{ij},\psi_i)$ for some open covering $\{ U_i\}$
of $U\times X$, where $\phi_{ij}\in \Ga(U_i\cap U_j,V^*(P))$,
$\psi_i\in\Ga(U_i,\gg^*(P))$ are such that
$-d^*\phi_{ij}=\psi_j-\psi_i$ over $U_i\cap U_j$,
$P$ is the universal $G$-bundle on $\MM$.
Similarly, every vector field $v\in T_{\MM}(U)$
can be represented by a Cech cocycle
$(\a_{ij},\nu_i)$, where $\a_{ij}\in\Ga(U_{ij},\gg(P))$,
$\nu_i\in\Ga(U_i, V(P))$ are such that $d\a_{ij}=\nu_j-\nu_i$.
In terms of these representatives the pairing between $T^*_{\MM}$
and $T_{\MM}$ takes form
\begin{equation}
\lan (\a_{ij},\nu_i),(\phi_{ij},\psi_i)\ran=
\Tr(\lan\a_{ij},\psi_j\ran+\lan\phi_{ij},\nu_i\ran),
\end{equation}
where $\Tr:H^1(U\times X,\O_{U\times X})\ra H^0(U,\O_U)$ is the morphism
induced by the trivialization of $\om_X$.

The map $H$ sends a 1-form $(\phi_{ij},\psi_i)$ to
the vector field represented by the cocycle $(t\circ d^*\phi_{ij},0)$
where $t$ is considered as a map $\gg^*\ra\gg$.
Since $d^*\phi_{ij}=\psi_i-\psi_j$ we have
$$H(\phi_{ij},\psi_i)=(0,d\circ t(\psi_i))\mod(\im(\de))$$
where $\de$ is the differential in the Cech complex of
$C^{\cdot}$. Note that since $d\circ t\circ d^*=0$, we
have $d\circ t(\psi_i)=d\circ t(\psi_j)$ over $U_i\cap U_j$,
hence we obtain the global section $d\circ t(\psi_{\cdot})\in
\Ga(U\times X, V(P))$.
It follows that
\begin{equation}
\lan H(\phi_{ij},\psi_i),(\phi'_{ij},\psi'_i)\ran=
\Tr(\lan\phi'_{ij},d\circ t(\psi_i)\ran)=
\Tr(\lan\psi'_i-\psi'_j, t(\psi_i)\ran.
\end{equation}

Let us consider the relative Atiyah extension for the universal bundle $P$:
$$0\ra\gg(P)\ra\AA(P)\ra p^*T_{\MM}\ra 0.$$
where $p:\MM\times X\ra\MM$ is the projection, $\AA(P)$ is the
bundle of relative  infinitesemal  symmetries  of $P$.
Then for sufficiently fine covering $\{ U_i\}$ a Cech cocycle
$(\a_{ij},\nu_i)$ representing a local vector field on
$\MM$  can  be   written   as   follows:   $\a_{ij}=D_j-D_i$,
$\nu_i=D_i(s)$ where $D_i\in\Ga(U_i,\AA(P))$, $s\in V(P)$
is the universal section, the symbol of $D_i$ is equal to the
restriction of a given vector field to $U_i$.
In particular, for a vector field
represented  by a cocycle $(0,\nu)$ where
$\nu\in\Ga(U\times X, V(P))$ we have $\nu=D(s)$ for some
$D\in\Ga(U\times X,\AA(P))$.

After these remarks we can start proving (\ref{Jacobi}).
Let us denote by $(\phi^h_{ij},\psi^h_i)$
Cech cocycles representing 1-forms $\om_h$ for $h=1,2,3$.
Let $D^h\in\Ga(U\times X,\AA(P))$ be infinitesemal
symmetries corresponding to $H(\om_h)$ so that
$D^h(s)=d\circ t(\psi_i)=t(\psi_i)(s)$ over $U_i$.
Then we have
$$H(\om_1)\cdot\lan H(\om_2),\om_3\ran=
D^1\cdot\Tr(\lan\psi^3_i-\psi^3_j, t(\psi^2_i)\ran)=
\Tr(\lan D^1(\psi^3_i-\psi^3_j), t(\psi^2_i)\ran+
\lan \psi^3_i-\psi^3_j, D^1(t(\psi^2_i))\ran).$$
On the other hand, it is easy to compute the commutator:
$$[H(\om_1),H(\om_2)]=(0,D^1D^2(s)-D^2D^1(s))=
(0,D^1(t(\psi^2_i)(s))-D^2(t(\psi^1_i)(s))).$$
Hence, we have
\begin{align*}
&\lan [H(\om_1),H(\om_2)], \om_3\ran=
\Tr\lan\phi^3_{ij},
D^1(t(\psi^2_i)(s))-D^2(t(\psi^1_i)(s))\ran=\\
&\Tr\lan\phi^3_{ij},D^1(t(\psi^2_i))(s)-D^2(t(\psi^1_i))(s)
+t(\psi^2_i)(D^1(s))-t(\psi^1_i)(D^2(s))\ran=
\\
&\Tr\lan\psi^3_i-\psi^3_j,D^1(t(\psi^2_i))-D^2(t(\psi^1_i))-
[t(\psi^1_i),t(\psi^2_i)]\ran.
\end{align*}
It follows that
\begin{align*}
&H(\om_1)\cdot\lan H(\om_2),\om_3\ran-
\lan [H(\om_1),H(\om_2)], \om_3\ran=\\
&\Tr(\lan D^1(\psi^3_i-\psi^3_j), t(\psi^2_i)\ran+
\lan\psi^3_i-\psi^3_j,t(D^2(\psi^1_i))+
[t(\psi^1_i),t(\psi^2_i)]\ran=\\
&
\Tr(\lan D^1(\psi^3_i-\psi^3_j), t(\psi^2_i)\ran+
\lan D^2(\psi^1_i), t(\psi^3_i-\psi^3_j)\ran+
\lan\psi^3_i-\psi^3_j,[t(\psi^1_i),t(\psi^2_i)]\ran)
\end{align*}
Since this is equal to
$D^1\cdot\Tr\lan\phi^3_{ij},dt(\psi^2_{\cdot})\ran-
\Tr\lan\phi^3_{ij},[D^1,D^2](s)\ran$
which is skew-symmetric in $i,j$,
we can also skew-symmetrize in $i,j$ the expression obtained above.
Then after adding terms obtained
by cyclic permutation of $1$, $2$ and $3$ we obtain the
trace of the following expression
$$\lan\psi^3_i-\psi^3_j,
[t(\psi^1_i),t(\psi^2_i)]+[t(\psi^1_j),t(\psi^2_j)]\ran+
cp(1,2,3).$$
Up to a coboundary this is equal to
$$\lan\psi^3_i,\psi^1_j,\psi^2_j\ran_t-
\lan\psi^3_j,\psi^1_i,\psi^2_i\ran_t+ cp(1,2,3),$$
where we denote $\lan x,y,z\ran_t=\lan x, [t(x),t(y)]\ran$.
From the fact that $t\in (S^2\gg)^*$ one can deduce
easily that $\lan \cdot,\cdot,\cdot\ran_t$ is $\gg$-invariant
and skew-symmetric. This implies the following identity:
$$\lan\psi^3_i,\psi^1_j,\psi^2_j\ran_t-
\lan\psi^3_j,\psi^1_i,\psi^2_i\ran_t+ cp(1,2,3)
=\lan\psi^1_i-\psi^1_j,\psi^2_i-\psi^2_j,\psi^3_i-
\psi^3_j\ran_t=
-\lan d^*\phi^1_{ij},d^*\phi^2_{ij},d^*\phi^3_{ij}\ran_t.$$
It  remains  to  notice  that  $\lan  d^*\phi^1,  d^*\phi^2,
d^*\phi^3\ran_t=0$
for any $\phi^1,\phi^2,\phi^3\in V^*(P)$.
Indeed, we have to show that for any triple of elements
$\varphi_1,\varphi_2,\varphi_3\in V^*$ and any $v\in V$ one has
$$\lan d^*_v\varphi_1, d^*_v\varphi_2, d^*_v\varphi_3\ran_t=0$$
where $d^*_v\varphi_h\in\gg^*$ is defined by
$d^*_v\varphi_h(x)=\varphi_h(x\cdot v)$, $h=1,2,3$.
Let $t=\sum_i x_i\otimes y_i$. Then
\begin{align*}
& \lan d^*_v\varphi_1,d^*_v\varphi_2,d^*_v\varphi_3\ran_t=
\lan d^*\varphi_1,
\sum_{i,j}[d^*\varphi_2(y_i)x_i,d^*\varphi_3(y_j)x_j]\ran=\\
&\sum_{i,j}\varphi_1([x_i,x_j]\cdot v)\varphi_2(y_i\cdot v)
\varphi_3(y_j\cdot v).
\end{align*}
Now
$$\sum_{i,j} [x_i,x_j]\cdot v\otimes y_i\cdot v\otimes
y_j\cdot v=
\sum_{i,j}x_ix_jv\otimes y_iv\otimes y_jv-
\sum_{i,j}x_jx_iv\otimes y_iv\otimes y_jv=0$$
since $\sum_i x_iv\otimes y_iv=0$.
\ed

Notice that the condition $t_*=0$ is usually not satisfied
when $\gg$ is simple. However, for example if $S^2(V)$
is irreducible and $\gg$ is simple then we necessarily
have $t_*=\la\cdot\id$ for some scalar $\la$. It follows
that we can replace $G$ by its product $G\times\G_m$ with one-dimensional
torus, $t$ by its sum with the appropriate multiple of the
square of the generator of $\Lie(\G_m)$, so that for the new
tensor $t'$ and the same representation $V$ (on which $\G_m$ acts
via identity character) the condition $t'_*=0$ will be satisfied.
In the case $G=\GL_{r_1}\times \GL_{r_2}$ and $V=\Mat(r_1,r_2)$
the tensor $t$ is equal to $(t_1, -t_2)$ where
$t_1=\sum  E_{ij}\otimes  E_{ji}$  is the standard symmetric
invariant tensor for $\gl_{r_1}$, $t_2$ is the similar
tensor for $\gl_{r_2}$.
An interesting case is $G=\GSp_{2r}$, the group of invertible
matrices preserving the symplectic from up to a scalar.
In   this   case   we  can  take  $V$  to  be  the  standard
representation  of  $G$  of  rank  $2r$,  then  $S^2(V)$  is
isomorphic to the adjoint representation of $\Sp_{2r}$
on which $\GSp_{2r}$ acts naturally, and one can easily
find a non-zero invariant tensor $t\in (S^2\gg)^{\gg}$ with $t_*=0$
(in fact, such $t$ is unique up to a constant).
The corresponding moduli stack is the stack of the following
data: a vector bundle $E$ together with a symplectic form
$$E\otimes E\ra L$$
inducing an isomorphism $E\simeq E^*\otimes L$, and a section
$s:\O_X\ra E$. In particular, taking $E$ to be a fixed bundle
with a symplectic form as above, we can consider sometimes
the appropriate quotient space of $\P H^0(E)$ by the group
of $\GSp$-automorphisms of $E$
as a Poisson substack in the above stack.
More precisely, we can define the stability condition for
such pairs $(E, s)$ depending on a parameter $\tau$:
the only difference with the case of $\GL$ is that one
should consider totally isotropic subbundles of $E$.
Then for $\tau=\mu(E)+\epsilon$ we have the Casimir morphism
from such moduli space to the stack of semistable
$\GSp$-bundles, hence its fibers inherit the Poisson structure.
For example, if $E_0$ is a stable bundle of degree 2
then there is a natural $\GSp_4$-structure on $E=E_0\oplus E_0$
such that both summands are totally isotropic. Then
the $\tau$-stability condition (with $\tau=\mu(E_0)+\epsilon)$
allows only sections $s\in H^0(E)=
H^0(E_0)\oplus H^0(E_0)$ with non-zero projections to both summands.
Hence, the space of such sections up to the action of symplectic
automorphisms of $E$ is $S^2\P H^0(E)$, so we get a Poisson
structure on the latter variety.

\end{document}